\documentclass[reqno, oneside]{amsart}
\usepackage{style}

\begin{document}
\title[Recurrence coefficients]{On recurrence coefficients of Steklov measures}
\author{R.~V.~Bessonov}

\address{St.Petersburg State University ({\normalfont \hbox{29b}, 14th Line V.O., 199178, St.Petersburg,  Russia}) and St.Petersburg Department of Steklov Mathematical Institute of Russian Academy of Science ({\normalfont 27, Fon\-tan\-ka, 191023, St.Petersburg, Russia})}
\email{bessonov@pdmi.ras.ru}

\thanks{The work is supported by RFBR grant mol\_a\_dk 16-31-60053 and by ``Native towns'', a social investment program of PJSC ``Gazprom Neft''}

\begin{abstract}
A measure $\mu$ on the unit circle $\T$ belongs to Steklov class $\Ss$ if its density $w$ with respect to the Lebesgue measure on $\T$ is strictly positive: $\inf_{\T} w > 0$. Let $\mu$, $\mu_{-1}$ be measures on the unit circle $\T$ with real recurrence coefficients $\{\alpha_k\}$, $\{-\alpha_k\}$, correspondingly.  If $\mu \in \Ss$ and $\mu_{-1} \in \Ss$, then partial sums $s_k=\alpha_0+ \ldots + \alpha_k$ satisfy the discrete Muckenhoupt condition  $\sup_{n > \ell\ge 0} \bigl(\frac{1}{n - \ell}\sum_{k=\ell}^{n-1} e^{2s_k}\bigr)\bigl(\frac{1}{n - \ell}\sum_{k=\ell}^{n-1} e^{-2s_k}\bigr) < \infty$. 
\end{abstract}

\maketitle

\section{Introduction}\label{s1}

Every probability measure $\mu$ on the unit circle $\T = \{z \in \C:\; |z| = 1\}$ of the complex plane $\C$ generates the family of monic orthogonal polynomials $\Phi_n$ satisfying the recurrence relations
$$
\Phi_{n+1}(z) = z \Phi_n(z) -\bar\alpha_{n} \Phi_n^\ast(z), \qquad z \in \T, \quad n \ge 0, \quad \Phi_0 \equiv 1,
$$
where $\Phi_n^\ast$ are the ``reversed'' polynomials defined by $\Phi_n^\ast(z) = z^n\ov{\Phi_n(1/\bar z)}$. The recurrence coefficients $\alpha_n=-\ov{\Phi_{n+1}(0)}$ are completely determined by the measure~$\mu$; in the non-trivial case where $\mu$ is supported on an infinite set, we have $|\alpha_n|<1$ for all $n \ge 0$. Any sequence of complex numbers $\alpha_n$ with $|\alpha_n|<1$ arises as the sequence of recurrence coefficients of a unique non-trivial probability measure on $\T$. A classical problem in the theory of orthogonal polynomials on the unit circle \cite{Simonbook} is to relate properties of probability measures $\mu$ to properties of their recurrence coefficients $\{\alpha_n\}$.
 
\medskip

In this paper we study recurrence coefficients of probability measures on $\T$ from Steklov class. Denote by $m$ the Lebesgue measure on $\T$ normalized by $m(\T) = 1$. A measure $\mu = w\,dm +\mu_s$ belongs to the Steklov class $\Ss$ if the density $w$ of its absolutely continuous part is strictly positive: 
\begin{equation*}
\inf_{z \in \T} w(z)>0.
\end{equation*} 
One version of famous Szeg\"o theorem says that 
\begin{equation}\label{eq8}
\prod_{k=0}^{\infty}(1-|\alpha_k|^2) = \exp\left(\int_{\T}\log w\,dm\right)
\end{equation}
for every probability measure $\mu = w\, dm + \mu_s$ on $\T$. If $\mu$ is a measure from Steklov class~$\Ss$, then $\log w \in L^1(\T)$, hence the product in the left hand side converges to a non-zero number and the recurrence coefficients of $\mu$ obey Szeg\"o condition $\sum|\alpha_k|^2 < \infty$. Another classical result, Baxter theorem, says that $\sum|\alpha_k| < \infty$ for the recurrence coefficients $\alpha_k$ of a probability measure $\mu$ on $\T$ if and only if $\mu$ is of the form $\mu = w\,dm$ for a strictly positive weight $w$ such that $\sum |\hat w(k)| <\infty$. Here and below $\hat w(k) = \int_{\T} w(z)\bar z^k \,dm(z)$, $k \in \Z$, denote the moments of $w$. See Chapters~2, 5 in \cite{Simonbook} for the proofs of Szeg\"o and Baxter theorems. Summarizing, condition $\sum|\alpha_k|^2 < \infty$ is necessary, while condition $\sum|\alpha_k| < \infty$ is sufficient for recurrence coefficients $\{\alpha_k\}$ to generate a measure $\mu\in \Ss$.

\medskip

Further information on recurrence coefficients of Steklov measures could be extracted from Rahmanov example solved the classical Steklov problem. The original question by Steklov asks if a sequence of orthogonal polynomials $P_n$ on the interval $[-1,1]$ generated  by a strictly positive weight $w$ on $[-1,1]$ is pointwise bounded: 
$$
\sup_n |P_n(x)| < \infty, \qquad x \in (-1,1).
$$ 
This question and closely related issues attracted a lot of attention, see detailed review \cite{Suetin}. The negative answer was given by Rahmanov \cite{Rahmanov} in 1979. After transferring the problem to the unit circle, he constructed a strictly positive weight $w$ on $\T$ such that $\sup_{n\ge 0}|\Phi_n(1)| = \infty$ for the orthogonal polynomials $\Phi_n$ it generates. This weight $w$ can be chosen to be symmetric \cite{Rahmanov2} with respect to the real line: $w(z) = w(\bar z)$ for almost all $z \in \T$. Note that for every symmetric weight $w$ its orthogonal polynomials $\Phi_n$ satisfy $\Phi_n(1) = \Phi_n^\ast(1)$, hence 
\begin{equation}\label{eq6}
\Phi_n(1) = \prod_{k = 0}^{n-1}(1-\alpha_k), \qquad n\ge 1.
\end{equation}
This formula implies $\inf_{n \ge 0} s_n = -\infty$ for the partial sums $s_n = \alpha_0 + \ldots +\alpha_n$ of recurrence coefficients $\{\alpha_k\}$ of the measure $w\,dm \in \Ss$ constructed in Rahmanov example. On the other hand, the Steklov bound 
\begin{equation}\label{eq5}
\max_{z \in \T}|\Phi_n(z)| \le \sum_{k = 0}^{n} |\hat\Phi_n(k)| \le \|\Phi_n\|_{L^2(\T)} \sqrt{n+1} \le (\inf\nolimits_\T w)^{-1}\|\Phi_n\|_{L^2(\mu)} \sqrt{n+1},  
\end{equation} 
yields the estimate  $s_n \ge -\frac{1}{2}\log n + c$ for all $n \ge 1$ and a constant $c$ independent of $n$ (to see this, use $\|\Phi_n\|_{L^2(\mu)} \le 1$, formula \eqref{eq6}, and the fact that $\sum|\alpha_k|^2 < \infty$). Recent advances in the area show that the Steklov bound is optimal in a natural sense \cite{ADT}, \cite{DK}. In particular, it follows from Theorem~4 in \cite{ADT} and formula \eqref{eq6} that for every positive sequence $\{\eps_k\}$ arbitrarily slowly tending to zero one can find a measure $\mu\in \Ss$ such that $s_{n_k} \le -\frac{1}{2}\log (\eps_k n_k)$ for some infinite increasing sequence $\{n_k\}$ of positive integers. See also \cite{Denisov2016}, \cite{DR17} for discussion of Rahmanov example and the corresponding recurrence coefficients.

\medskip

In this paper we develop a method allowing to control oscillations of the sequence of partial sums, $\{s_n\}$, of recurrence coefficients of Steklov measures. Our main result is the following theorem.
   
\begin{Thm}\label{t1}
Let $\mu$, $\mu_{-1}$ be measures on the unit circle $\T$ with real recurrence coefficients $\{\alpha_k\}$, $\{-\alpha_k\}$, correspondingly.  If $\mu \in \Ss$ and $\mu_{-1} \in \Ss$, then
\begin{equation}\label{eq55}
\sup_{n > \ell \ge 0} \, \left(\frac{1}{n - \ell}\sum\nolimits_{k=\ell}^{n-1} e^{2s_k}\right) \left(\frac{1}{n - \ell}\sum\nolimits_{k=\ell}^{n-1} e^{-2s_k}\right) < \infty,
\end{equation}
where $s_k = \alpha_0 + \ldots + \alpha_k$ for integer $k \ge 0$.   
\end{Thm}
The Muckenhoupt class $A_2(\R)$ on the real line $\R$ consists of functions $g$ such that 
\begin{equation}\label{eq66}
\sup_{I \subset\R} \left(\frac{1}{|I|}\int_{I} g(t)\,dt\right)\left(\frac{1}{|I|}\int_{I} \frac{1}{g(t)}\,dt\right) < \infty,
\end{equation}
where the supremum is taken over all intervals $I \subset \R$. A similarity between relations~\eqref{eq55} and \eqref{eq66} explains the name 
``discrete Muckenhoupt condition'' we use for referring to~\eqref{eq55}. One may observe that in the setting of the Baxter theorem the partial sums $s_k$ are uniformly bounded and hence relation \eqref{eq55} is obviously satisfied. Jensen's inequality implies that the sequence $s=\{s_k\}$ in Theorem \ref{t1} has bounded mean oscillation:
\begin{equation}\label{eq7}
\sup_{n>\ell\ge0}\frac{1}{n - \ell}\sum\nolimits_{k=\ell}^{n-1} \bigl|s_k - \langle s \rangle_{n,\ell}\bigr| < \infty, \qquad \langle s \rangle_{n,\ell} = \frac{1}{n - \ell}\sum\nolimits_{k=\ell}^{n-1}s_k. 
\end{equation}
According to John-Nirenberg inequality, sequences of bounded mean oscillation grow at most logarithmically. This agrees well with the Steklov bound \eqref{eq5}. 

\medskip
 
The measure $\mu_{-1}$ in Theorem~\ref{t1} is the orthogonality measure for the second kind polynomials generated by $\mu$. Given $\mu$, it is possible to construct the measure $\mu_{-1}$ not knowing the recurrence coefficients $\{\alpha_k\}$. Theorem~\ref{t1} then can be reformulated without referring to $\mu_{-1}$. For more details, see Section \ref{s3}.  

\medskip

The author wishes to thank Stanislav Kupin from University Bordeaux 1 who advised me to search for an analogue of Theorem 1 from \cite{B17} in the theory of orthogonal polynomials, inspiring this work.   

\section{Proof of Theorem \ref{t1}}\label{s2}
Let $\mu$ be a probability measure on the unit circle $\T$ supported on a set of infinitely many points, and let $\{\Phi_n\}_{n\ge 0}$ be the sequence of monic polynomials orthogonal with respect to $\mu$. Recall that the polynomials $\Phi_n$ are determined by relations 
$$
\deg\Phi_n = n, \qquad (\Phi_n, \Phi_k)_{L^2(\mu)} = 0 \mbox{ for }k\neq n, \qquad \hat\Phi_n(n) = 1,
$$
and could be obtained via Gram-Schmidt orthogonalization of $\{z^n\}_{n \ge 0}$. These polynomials satisfy the system of recurrence relations 
\begin{equation}\label{eq9}
\begin{cases}
\Phi_{n+1}(z) = z \Phi_n(z) -\bar\alpha_n \Phi_n^\ast(z), \quad &\Phi_0 \equiv 1,\\
\Phi_{n+1}^\ast(z) = -\alpha_n z \Phi_n(z) +   \Phi_n^\ast(z), \quad &\Phi^\ast_0 \equiv 1. 
\end{cases}
\end{equation}
The numbers $\alpha_n$, $n \ge 0$, are called the recurrence (or Schur/Verblunski/reflection) coefficients of the measure $\mu$. 
By definition, we have $\alpha_n = -\ov{\Phi_{n+1}(0)}$. For basic theory of orthogonal polynomials on the unit circle we refer the reader to book~\cite{Simonbook}. 

\medskip

Fix a probability measure $\mu$ on the unit circle $\T$ having real recurrence coefficients~$\{\alpha_k\}$. Let $\Phi_n$ be the monic orthogonal polynomials with respect to $\mu$. For $\alpha \in \R$ define 
$$
T(\alpha,z) = 
\begin{pmatrix}
z  &-\alpha \\
-\alpha z &1
\end{pmatrix},
\qquad 
T(\alpha) = 
\begin{pmatrix}
1 &-\alpha \\
-\alpha &1
\end{pmatrix},
\qquad
Q(\alpha) = 
\begin{pmatrix}
1 & 0 \\
-\alpha & 0
\end{pmatrix}.
$$ 
Then relations \eqref{eq9} yield
\begin{equation}\label{eq12}
\begin{pmatrix}\Phi_{n+1}(z) \\ \Phi_{n+1}^\ast(z)\end{pmatrix}= 
T(\alpha_{n},z)\cdot \ldots \cdot T(\alpha_0,z) \begin{pmatrix}1\\1\end{pmatrix}, 
\qquad n\ge0. 
\end{equation}
In particular, we have
$$
\Phi_{n+1}(1) = \Phi_{n+1}^\ast(1) = \prod_{k = 0}^{n}(1-\alpha_k), \quad n \ge 0.
$$
Below in Lemma \ref{l2} we present a formula in terms of $\{\alpha_k\}$ for derivatives $\Phi_{n+1}^{(j)}(1)$ of order $1 \le  j \le n+1$ evaluated at the point $1$. This formula will play a central role in our considerations. For a multi-index $\gamma=(\gamma_0,\ldots,\gamma_{n})$ of length $n+1$ with components $0$ and $1$, put
$$
\Pi(\gamma) =  \Pi_{\gamma_{n}}(\alpha_{n})\cdot \ldots \cdot\Pi_{\gamma_{1}}(\alpha_1)\Pi_{\gamma_{0}}(\alpha_0),
$$
where $\Pi_{0}(\alpha) = T(\alpha)$ and $\Pi_{1}(\alpha) = Q(\alpha)$ for all $\alpha \in \R$. Denote by $\M_{n,j}$ the set of all multi-indexes $\gamma=(\gamma_0,\ldots\gamma_{n})$ with components~$0, 1$ such that $\gamma_0 + \ldots + \gamma_n = j$. We start with a simple lemma.	
\begin{Lem}\label{l1} 
For all integers $n,j$ such that $0 \le j \le n+1$ we have
\begin{equation}\label{eq1}
\left(\begin{smallmatrix}\Phi_{n+1}^{(j)}(1) \\ {\Phi_{n+1}^\ast}^{\!(j)}\!(1)\end{smallmatrix}\right) 
=j!\sum_{\gamma \in \M_{n,j}}  \Pi(\gamma)\onne. 
\end{equation} 
\end{Lem}
\beginpf Formula \eqref{eq1} for $j = 0$ and all integers $n \ge 0$ is just relation \eqref{eq12} for $z=1$. For $j\ge 1$, we can differentiate the expression in formula \eqref{eq12} $j$ times and obtain
$$
\left(\begin{smallmatrix}\Phi_{n+1}^{(j)}(1) \\ {\Phi_{n+1}^{\ast}}^{\hspace{-3.5mm}(j)}\!(1)\end{smallmatrix}\right) 
= j!\sum_{\gamma \in \M_{n,j}} \partial^{\gamma}T(\alpha_n,z_n)\cdot \ldots \cdot T(\alpha_{0},z_{0})|_{\mathbf{z}},
\quad \partial^{\gamma} = \frac{\partial^{\gamma_{n}}}{\partial{z_{n}^{\gamma_{n}}}} \ldots \frac{\partial^{\gamma_0}}{\partial{z_{0}^{\gamma_0}}},
$$
where $\mathbf z = (1,1,\ldots, 1)$ is the vector in $\R^{n+1}$. After substitution $T(\alpha,z)|_{z=1} = T(\alpha)$ and $T(\alpha,z)' \equiv Q(\alpha)$, we get the desired proposition. \qed

\medskip

\noindent For integers $k_1, k_2$ and  a sequence of real numbers $\{c_k\}$ we denote 
$$
\psum_{k=k_1}^{k_2} c_k = \begin{cases}\sum_{k=k_1}^{k_2-1} c_k, &\mbox{ if }k_1<k_2, \\  0, &\mbox{ if }k_1 \ge k_2.\end{cases}
$$ 
It will be convenient to put $\alpha_{-1} = 0$ and to define the function (sequence) on integers,
$$
h: n \mapsto \prod_{k=-1}^{n-1} \frac{1 - \alpha_k}{1+\alpha_k}, \quad n \ge 0.
$$    
\begin{Lem}\label{l2}
For all integers $n,j$ such that $1 \le j \le n+1$, we have
\begin{equation}\label{eq2}
\Phi_{n+1}^{(j)}(1) = \frac{j!}{2^j}\prod_{t=0}^{n}(1-\alpha_t)\sum_{t_1 = 0}^{n}\psum_{t_2 = 0}^{t_1}\ldots 
\psum_{t_{j} = 0}^{t_{j-1}}\prod_{s=1}^{j}\left(1 + \frac{h(t_{s})}{h(t_{s-1})}\right),
\end{equation} 
where we put $t_{0} = n+1$.
\end{Lem}
\beginpf Take a multi-index $\gamma \in \M_{n,j}$. Let $0\le t_j< \ldots < t_1 \le n$ be the indexes such that $\gamma_{t_s} = 1$ for $1\le s \le j$. As in the statement of the Lemma, put $t_{0} = n+1$. Using identities
\begin{align*}
&T(\alpha)\onne = (1-\alpha)\onne, \quad T(\alpha)\onnem = (1+\alpha)\onnem,\\
&Q(\alpha)\onne = Q(\alpha)\onnem = \tfrac{1-\alpha}{2}\onne + \tfrac{1+\alpha}{2}\onnem,
\end{align*}
we obtain
\begin{align*}
\bigl\langle\Pi(\gamma)\onne \onnez\bigr\rangle 
&= \frac{1}{2^j}\prod_{s=1}^{j}\left(\prod\nolimits_{t=t_{s}}^{t_{s-1}-1}(1-\alpha_t) + \prod\nolimits_{t=t_{s}}^{t_{s-1}-1}(1+\alpha_t)\right)\cdot \prod_{t=-1}^{t_j-1}(1-\alpha_t),\\
&= \frac{1}{2^j}\prod_{t=0}^{n}(1-\alpha_t)\cdot \prod_{s=1}^{j}\left(1 + \frac{h(t_{s})}{h(t_{s-1})}\right),
\end{align*}
Summing up over all multi-indexes $\gamma\in \M_{n,j}$ and using Lemma \ref{l1}, we obtain formula \eqref{eq2}. \qed

\medskip

Denote by $H^2(\mu,n)$ the subspace in $L^2(\mu)$ consisting of all analytic polynomials of degree at most $n$.  Let $k_{\zeta,\mu,n}$ be the reproducing kernel in $H^2(\mu,n)$ at $\zeta \in \C$. Define $\pi_r = \prod_{k=-1}^{r-1}(1 - \alpha_k^2)$ for $r \ge 0$. The following well-known relation follows from the fact that the family $\bigl\{\tfrac{1}{\sqrt{\pi_r}} \Phi_r\bigr\}_{0 \le r \le n}$ is the orthonormal basis in $H^2(\mu,n)$:
\begin{equation}\label{eq3}
k_{\zeta,\mu,n}(z) = \sum_{r=0}^{n}\frac{1}{\pi_{r}} \Phi_r(z)\ov{\Phi_r(\zeta)}, \qquad z \in \C. 
\end{equation}
See Section 2.2 in \cite{Simonbook} for more details. For an integer $0 \le j \le n$ we denote by $\bar\partial^{j}k_{\xi,\mu,n}$ the derivative of order $j$ of the anti-analytic mapping $\zeta \mapsto k_{\zeta,\mu,n}$ evaluated at a point $\xi \in \C$.

\begin{Lem}\label{l7}
Let $\mu$ be a measure from Steklov class. Then there exists a constant $c$ such that $\|\bar\partial^{j}k_{1,\mu,n}\|_{L^2(\mu)}^{2} \le c \|\bar\partial^{j}k_{1,m,n}\|_{L^2(m)}^{2}$ for all integers $n \ge 0$ and  $0 \le j \le n$.   
\end{Lem}
\beginpf For every $n\ge 0$ consider the operator $T_{\mu,n}$ on $H^2(m,n)$ defined by 
$$
(T_{\mu,n}f,g)_{L^2(m)} = \int_{\T}f\bar g\,d\mu.
$$
Since the measure $\mu = w\,dm + \mu_s$ belongs to the Steklov class $\Ss$, we have
\begin{equation}\label{eq21}
(T_{\mu,n}f,f)_{L^2(m)} \ge \int_{\T}|f|^2w\,dm \ge \inf_{z \in \T}w(z) \cdot \|f\|^{2}_{L^2(\T)}.
\end{equation}
In particular, the operators $T_{\mu,n}$, $n \ge 0$, are invertible and the supremum 
$$
c = \sup_{n\ge 1}\|T_{\mu,n}^{-1}\|_{H^2(m,n) \to H^2(m,n)}
$$ 
is finite. Take a point $\zeta \in \C$ and consider $T_{\mu,n}^{-1} k_{\zeta,m,n}$ as an element of $H^2(\mu,n)$. For every polynomial $f \in H^2(\mu,n)$ we have 
$$
(f, T_{\mu,n}^{-1} k_{\zeta,m,n})_{H^2(\mu,n)} = (T_{\mu,n} f, T_{\mu,n}^{-1} k_{\zeta,m,n})_{L^2(m)} = f(\zeta) = (f, k_{\zeta,\mu,n})_{H^2(\mu,n)}. 
$$
It follows that $k_{\zeta,\mu,n} = T_{\mu,n}^{-1} k_{\zeta,m,n}$. Differentiating this relation with respect to $\bar \zeta$ at $\zeta=1$, we obtain 
$$
\bar\partial^j k_{1,\mu,n} = T_{\mu,n}^{-1} \bar\partial^j k_{1,m,n}, \qquad 0 \le j \le n.
$$
From here and the definition of $c$ we see that $\|\bar\partial^{j}k_{1,\mu,n}\|_{L^2(\mu)}^{2} \le c \|\bar\partial^{j}k_{1,m,n}\|_{L^2(m)}^{2}$ for all~$n$ and $j$, as required. \qed

\begin{Lem}\label{l3}
For all integers $n \ge j \ge 1$ we have 
\begin{equation}\label{eq4}
\|\bar\partial^{j}k_{1,\mu,n}\|_{L^2(\mu)}^{2} 
= \sum_{r=j}^{n}h(r)\left(\frac{j!}{2^j}\psum_{t_1 = 0}^{r}\psum_{t_2 = 0}^{t_1}\ldots 
\psum_{t_j = 0}^{t_{j-1}}\prod_{s=1}^{j}\left(1+\frac{h(t_{s})}{h(t_{s-1})}\right)\right)^{2},
\end{equation}
and $\|k_{1,\mu,n}\|_{L^2(\mu)}^{2} = \sum_{r=0}^{n}h(r)$ for the case where $j=0$. 
\end{Lem}
\beginpf For $j = 0$ we have 
$$
\|k_{1,\mu,n}\|_{L^2(\mu)}^{2} = \sum_{r=0}^{n}\frac{1}{\pi_{r}} |\Phi_{r}(1)|^{2} = \sum_{r=0}^{n}\frac{1}{\pi_{r}}\prod_{k=-1}^{r-1}(1-\alpha_k)^{2} = \sum_{r=0}^{n}h(r).
$$
In the case where $j \ge 1$ we differentiate \eqref{eq3} and obtain 
$$
\|\bar\partial^{j}k_{1,\mu,n}\|_{L^2(\mu)}^{2} = \bar\partial^{j} k_{1,\mu,n}^{(j)}(1) = \sum_{r=j}^{n} \frac{1}{\pi_{r}}|\Phi_{r}^{(j)}(1)|^2. 
$$
Formula \eqref{eq4} now follows readily from formula \eqref{eq2}. \qed

\medskip

\noindent For integers $r \ge 0$, $j\ge1$ we define 
$$
\kappa_j(r) = j! \psum_{t_1 = 0}^{r}\psum_{t_2 = 0}^{t_1}\ldots 
\psum_{t_{j} = 0}^{t_{j-1}} 1 = 
\begin{cases}
r!/(r-j)! & \mbox{if } r\ge j,\\
0 & \mbox{if } r< j.
\end{cases}
$$
\begin{Lem}\label{l5}
Let $n$, $j$ be integer numbers such that $1 \le j \le n/2$. 
Denote by $n_j$ the integer part of the number $(1 - \frac{1}{j+1})n$. We have 
$\frac{\kappa_j(n)}{\kappa_j(n_j)} \le C$
for a universal constant $C$ independent of $n$ and $j$.
\end{Lem}
\beginpf By Stirling formula, the fraction $\frac{\kappa_j(n)}{\kappa_j(n_j)} = \frac{n!(n_j-j)!}{n_j!(n-j)!}$ is comparable to 
\begin{align*}
\frac{n^n (n_j - j)^{n_j - j}\sqrt{n(n_j - j)}}{{n_j}^{n_j}(n - j)^{n - j}\sqrt{n_j(n - j)}}. 
\end{align*}
We can assume that $n \ge 10$. Then for all $1 \le j \le n/2$ we have 
$$
\frac{n}{n_j} \le \frac{n}{(1- \frac{1}{j+1})n -1} \le \frac{n}{n/2-1} \le 3, \qquad \frac{n_j - j}{n - j} \le 1. 
$$ 
So it suffices to show that the quantity 
\begin{equation*} 
A_{n,j} = n \log n + (n_j-j)\log(n_j-j) -n_j\log n_j - (n-j)\log(n-j)  
\end{equation*}
is bounded from above by a constant do not depending on $n \ge 10$ and $1 \le j \le n/2$. For such indexes $n,j$ we have
$$
\frac{n-n_j}{n-j} \le \frac{n+1}{(j+1)(n-j)} \le \min\left(\frac{n+1}{2(n-1)}, \frac{n+1}{3(n-n/2)}\right) \le \frac{3}{4}. 
$$
Let $c$  be a constant such that $|\log (1 + x)| \le c|x|$ for all $|x|\le\frac{3}{4}$. Then
\begin{align*}
A_{n,j}
&=-n_j \log\left(1 + \frac{j(n - n_j)}{n(n_j - j)}\right) - j\log \left(1 - \frac{n - n_j}{n-j}\right) - (n-n_j)\log(1 - \tfrac{j}{n})\\
&\le c\frac{j(n-n_j)}{n-j} + c\frac{j(n-n_j)}{n} \le 2c \frac{j(n-n_j)}{n-j} \le 2c \frac{j(n+1)}{(j+1)(n-j)} \le 8c.
\end{align*}
The lemma follows. \qed

\medskip

\noindent{\bf Proof of Theorem \ref{t1}.} We will prove that for every pair of integers $n>\ell\ge0$ one can find an integer $j \le \frac{n}{2}$ depending on $\ell$ such that 
$$
\left(\frac{1}{n - \ell}\sum_{k=\ell}^{n-1} e^{2s_k}\right)\left(\frac{1}{n - \ell}\sum_{k=\ell}^{n-1} e^{-2s_k}\right) \le 
c\frac{\|\bar\partial^{j}k_{1,\mu,n}\|_{L^2(\mu)}^{2}\|\bar\partial^{j}k_{1,\mu_{-1},n}\|_{L^2(\mu_{-1})}^{2}}{\|\bar\partial^{j}k_{1,m,n}\|_{L^2(m)}^{4}},
$$
where $c$ is a universal constant. Since the right hand side is uniformly bounded in $n$, $j$ by Lemma \ref{l7}, this is sufficient for the proof of the statement. 

\medskip

Let $n$, $j$ be integer numbers such that $n \ge 10$ and $1 \le j \le n/2$. 
Denote by $n_j$ the integer part of the number $(1 - \frac{1}{j+1})n$. 
Using Lemma \ref{l3} and Jensen's inequality, we obtain
\begin{align*}
\|\bar\partial_{1}^{j}k_{1,\mu,n}\|_{L^2(\mu)}^{2} 
&= \sum_{r=j}^{n}h(r)\left(\frac{j!}{2^j}\psum_{t_1 = 0}^{r}\psum_{t_2 = 0}^{t_1}\ldots 
\psum_{t_j = 0}^{t_{j-1}}\prod_{s=1}^{j}\left(1+\frac{h(t_{s})}{h(t_{s-1})}\right)\right)^{2},\\
&\ge \sum_{r=n_j}^{n}h(r)\left(\frac{j!}{2^j}\psum_{t_1 = 0}^{n_j}\psum_{t_2 = 0}^{t_1}\ldots 
\psum_{t_j = 0}^{t_{j-1}}\prod_{s=1}^{j}\left(1+\frac{h(t_{s})}{h(t_{s-1})}\right)\right)^{2},\\
&\ge\frac{e^{G_{h,n,j}}\kappa_j(n_j)^2}{4^j} \sum_{r=n_j}^{n}h(r),
\end{align*}
where
$$
G_{h,n,j} = \frac{2}{\kappa_j(n_j)}\psum_{t_1 = 0}^{n_j}\psum_{t_2 = 0}^{t_1}\ldots 
\psum_{t_j = 0}^{t_{j-1}}\sum_{s=1}^{j}\log\left(1+\frac{h(t_{s})}{h(t_{s-1})}\right).
$$
The same consideration applies to the triple $\mu_{-1}$, $\{-\alpha_k\}$, $\frac{1}{h}$, and gives
$$
\|\bar\partial_{1}^{j}k_{1,\mu_{-1},n}\|_{L^2(\mu_{-1})}^{2} \ge \frac{e^{G_{1/h,n,j}}\kappa_j(n_j)^2}{4^j} \sum_{r=n_j}^{n}\frac{1}{h(r)}.
$$
Since $\log(1+a) + \log(1 + a^{-1}) \ge \log 4$ for every $a>0$,
we see that
$$
\frac{e^{G_{h,n,j}}}{4^j} \cdot \frac{e^{G_{1/h,n,j}}}{4^j} \ge \frac{e^{2 \sum_{s = 1}^{j}\log 4}}{16^j} = 1.
$$
It follows that
$$
\|\bar\partial_{1}^{j}k_{1,\mu,n}\|_{L^2(\mu)}^{2}\|\bar\partial_{1}^{j}k_{1,\mu_{-1},n}\|_{L^2(\mu_{-1})}^{2} \ge 
\kappa_j(n_j)^4\left(\sum_{r=n_j}^{n}h(r)\right)\left(
\sum_{r=n_j}^{n}\frac{1}{h(r)}\right).
$$
On the other hand,
\begin{align*}
\|\bar\partial_{1}^{j}k_{1,m,n}\|_{L^2(m)}^{2} 
&= \sum_{r=j}^{n}\left(j!\psum_{t_1 = 0}^{r}\psum_{t_2 = 0}^{t_1}\ldots 
\psum_{t_j = 0}^{t_{j-1}}1 \right)^{2} = \sum_{r=j}^{n} \left(\frac{r!}{(r - j)!}\right)^2 \\
&= \left(\frac{n!}{(n - j)!}\right)^2 \left(1 + \sum_{r=j}^{n-1} \left(\frac{r-j+1}{r+1}\cdot \ldots \cdot \frac{n-j}{n}\right)^2\right) \\
&\le \left(\frac{n!}{(n - j)!}\right)^2 \left(1 + \sum_{r=j}^{n-1} \left(\frac{n-j}{n}\right)^{2(n-r)}\right)\\
&\le \left(\frac{n!}{(n - j)!}\right)^2 \sum_{s=0}^{\infty} \left(\frac{n-j}{n}\right)^{2s} \le \frac{n}{j}\left(\frac{n!}{(n - j)!}\right)^2. 
\end{align*}
From here we see that 
\begin{align*}
\frac{\|\bar\partial_{1}^{j}k_{1,\mu,n}\|_{L^2(\mu)}^{2}\|\bar\partial_{1}^{j}k_{1,\mu_{-1},n}\|_{L^2(\mu_{-1})}^{2}}{\|\partial_{1}^{j}k_{1,m,n}\|_{L^2(m)}^{4}}
&\ge \frac{j^2 \kappa_j(n_j)^4}{n^2 \kappa_j(n)^4} \left(\sum_{r=n_j}^{n}h(r)\right)\left(\sum_{r=n_j}^{n}\frac{1}{h(r)}\right)\\
&\ge \frac{1}{4C^{4}}\left(\frac{1}{n-n_j}\sum_{r=n_j}^{n}h(r)\right)\left(\frac{1}{n-n_j}\sum_{r=n_j}^{n}\frac{1}{h(r)}\right),
\end{align*}
where we used Lemma \ref{l5} and the fact that 
$$
\frac{1}{n-n_j} \le \frac{1}{n - (1 - \frac{1}{j+1})n} = \frac{j+1}{n} \le \frac{2j}{n}. 
$$
Now Lemma \ref{l7} applied to measures $\mu$, $\mu_{-1}$ from $\Ss$ yields the inequality
\begin{equation}\label{eq22}
\left(\frac{1}{n-\ell}\sum\nolimits_{k=\ell}^{n}h(k)\right)\left(\frac{1}{n-\ell}\sum\nolimits_{k=\ell}^{n}\frac{1}{h(k)}\right)\le 4C^4
\end{equation}
for all pairs $n, \ell$ such that $\ell=[(1 - \frac{1}{j+1})n]$ for some $1 \le j \le \frac{n}{2}$, where $[a]$ denotes the integer part of a real number $a>0$. By Lemma \ref{l7} and Lemma \ref{l3} for $\mu$, $\mu_{-1}$, and $j = 0$, inequality \eqref{eq22} holds in the case where $\ell=0$ as well. Now take arbitrary integers $n>\ell\ge0$ and find maximal integer $j \le \frac{n}{2}$ such that $\ell^* = [(1 - \frac{1}{j+1})n] \le \ell$. By construction, $n-\ell$ is comparable to $n-\ell^*$ with absolute constants. Hence, we can estimate
\begin{align*}
\frac{1}{(n-\ell)^2}\left(\sum_{k=\ell}^{n}h(k)\right)\left(\sum_{k=\ell}^{n}\frac{1}{h(k)}\right)	
\le
\frac{C_1}{(n-\ell^*)^2}\left(\sum_{k=\ell^*}^{n}h(k)\right)\left(\sum_{k=\ell^*}^{n}\frac{1}{h(k)}\right)
\le C_2.
\end{align*} 
It follows that inequality \eqref{eq22} holds for all $n$, $\ell$ and some new constant $C$. 
By Szeg\"o theorem~\eqref{eq8} and the definition of Steklov class $\Ss$, we have $\sum_{k \ge 0}|\alpha_k|^2 < \infty$. In particular, $\sup_k|\alpha_k| < 1$ and 
$$
|2s_n+\log h(n)| \le \sum_{k = 0}^{n} \left|2\alpha_k - \log \frac{1+\alpha_k}{1-\alpha_k}\right| \le c\sum_{0}^{\infty}|\alpha_k|^2,
$$ 
for a constant $c$ do not depending on $n$. Now the discrete Muckenhoupt condition~\eqref{eq55}  follows from formula \eqref{eq22}. \qed

\medskip

\section{Reformulation of Theorem \ref{t1}. Negative recurrence coefficients}\label{s3}
Theorem~\ref{t1} could be reformulated in a way avoiding the usage of the ``second kind'' measure $\mu_{-1}$. For this we need the definition of the Hilbert transform $Hf$ of a function $f \in L^1(\T)$:
$$
Hf(z) = \dashint \frac{f(\xi)}{1-\bar\xi z}\,dm(\xi), \qquad z \in \T.
$$
It is known that the principal value integral in formula above converges for almost all $z\in \T$, see Section III.1 in \cite{Garnett}. 
\addtocounter{Thm}{-1}
\renewcommand{\theThm}{\arabic{Thm}$'$} 
\begin{Thm}\label{t1p}
Let $w\ge 0$ be a function on $\T$ such that $w(z) = w(\bar z)$ for almost all~$z$. Assume that $\inf_{\T}w>0$ and $\sup_{\T} w < \infty$. If, moreover, $Hw$ is bounded on $\T$, then 
$$
\sup_{n > \ell \ge 0} \, \left(\frac{1}{n - \ell}\sum\nolimits_{k=\ell}^{n-1} e^{2s_k}\right) \left(\frac{1}{n - \ell}\sum\nolimits_{k=\ell}^{n-1} e^{-2s_k}\right) < \infty,
$$
where $s_k = \alpha_0 + \ldots + \alpha_k$, $k \ge 0$, denote the partial sums of recurrence coefficients of $\mu = w\,dm$.   
\end{Thm}
\renewcommand{\theThm}{\arabic{Thm}}
\noindent Both Szeg\"o and Baxter conditions,
$$
\sum\limits_{k\ge 0}|\alpha_k|^2 < \infty, \qquad  \sum\limits_{k\ge 0}|\alpha_k| < \infty,
$$
are invariant under the multiplication of the sequence $\{\alpha_k\}$ by a number $\lambda$ of unit modulus. The  situation for the Steklov class is completely different.
\begin{Prop}\label{p1}
Let $\mu$, $\mu_{-1}$ be probability measures on the unit circle $\T$ with recurrence coefficients $\{\alpha_k\}$, $\{-\alpha_k\}$, correspondingly. The following are equivalent:
\begin{itemize}
\item[(a)] $\mu \in \Ss$ and $\mu_{-1} \in \Ss$, 
\item[(b)] $\mu = w\,dm$ for a weight $w$ on $\T$ such that $\inf_{\T}w>0$, $\sup_{\T} w < \infty$, and the Hilbert transform $Hw$ is bounded on $\T$. 
\end{itemize}
\end{Prop} 
\noindent It is clear from Proposition~\ref{p1} that Theorem~\ref{t1p} is equivalent to Theorem \ref{t1}. The only thing to note here is that we have $w(z) = w(\bar z)$ for a strictly positive weight $w$ on $\T$ and almost all $z \in \T$ if and only if the recurrence coefficients of the measure $\mu = w\, dm$ are real. The latter follows from the fact that $w$ could be weakly approximated by the sequence of weights $w_n = \frac{1}{|\Phi_n^*|^2}$ (see Theorem 1.7.8 in \cite{Simonbook}), which satisfy relations $w_n(z) = w_n(\bar z)$, $z \in \T$.  

\medskip

In the proof of Proposition \ref{p1} we will use a couple of facts from the theory of harmonic functions. First, a measure $\mu$ on $\T$ has bounded density with respect to the Lebesgue measure $m$ on $\T$ if and only if the harmonic extension of $\mu$ to the open unit disk $\D = \{z \in \C:\; |z|<1\}$,
$$
\mathcal{P\mu}(z) = \int_{\T}\frac{1-|z|^2}{|1-\bar\xi z|^2}\,d\mu(\xi), \qquad z \in \D,
$$ 
is bounded. Second, let $w \in L^1(m)$ and denote by $F$ the analytic function in $\D$ such that $\mathcal P\mu(z) = \Re F(z)$, $z \in \D$, where $\mu = w\,dm$. Then $\Im F$ is bounded in $\D$ if and only if the Hilbert transform of $w$ is bounded on $\T$. For the proof of these facts, see, e.g., Sections I.3 and III.1 in \cite{Garnett}. 

\medskip

\noindent{\bf Proof of Proposition \ref{p1}.} Consider the analytic function $F$ in the open unit disk such that its real part $\Re F$ coincides with the harmonic extension of $\mu$,
\begin{equation}\label{eq11}
\Re F(z) = \int_{\T}\frac{1-|z|^2}{|1-\bar\xi z|^2}\,d\mu(\xi), \qquad |z|<1,
\end{equation} 
and $\Im F(0) = 0$. The analytic function $\frac{1}{F}$ has positive real part and equals one at the origin, hence there exists a probability measure $\nu$ on $\T$ such that $\Re\frac{1}{F}$ is the harmonic extension of $\nu$. By Theorem 3.2.14 in \cite{Simonbook}, the recurrence coefficients of $\nu$ coincide with $\{-\alpha_k\}$, that is, $\nu = \mu_{-1}$. Thus, the harmonic extension of $\mu_{-1}$ to the open unit disk 
has the form
\begin{equation}\label{eq10}
\Re\left(\frac{1}{F(z)}\right) = \frac{\Re F(z)}{(\Re F(z))^2 + (\Im F(z))^2}, \qquad |z| < 1.
\end{equation}
Consider the case where $\mu$, $\mu_{-1}$ belong to the Steklov class. Then $\Re\frac{1}{F}$, the harmonic extension of $\mu_{-1}$, is such that $\inf_{|z|<1}\Re\frac{1}{F(z)} > 0$, and we see from \eqref{eq10} that $\Re F$ must be bounded in $\D$. We also have $\inf_{|z|<1}\Re F(z) > 0$ by \eqref{eq11}, hence the measure $\mu$ is absolutely continuous, $\mu = w\,dm$, and $\inf_{\T}w>0$, $\sup_{\T} w < \infty$. Moreover, since  both functions $\Re F$, $\Re\frac{1}{F}$ are separated from zero, the function $\Im F$ is bounded in $\D$, see \eqref{eq10}. It follows that the Hilbert transform of $w$ is bounded on the unit circle $\T$.   

Now consider a measure $\mu$ as in $(b)$. We obviously have $\mu \in \Ss$. Let $F$ be the analytic function defined by \eqref{eq11}. Then  $\Re F$ is  bounded in $\D$ and we have $\inf_{z \in \D} \Re F(z) > 0$. Since the Hilbert transform of $w$ is bounded on $\T$, the function $\Im F$ is bounded in $\D$. Hence $\Re\frac{1}{F}$ is a strictly positive bounded function in $\D$. It follows that $\mu_{-1}$ is an absolutely continuous measure on $\T$ whose density is bounded and separated from zero, in particular, we have $\mu_{-1} \in \Ss$. \qed

\medskip

\noindent{\bf Remark.} It is not known to the author if there exists a measure $\mu \in \Ss$ with real recurrence coefficients $\{\alpha_k\}$ that do not obey the discrete Muckenhoupt condition~\eqref{eq55}. In a similar situation \cite{B17} on the real line the Muckenhoupt condition  holds for all measures of the form $\mu = w\,dx$, where $w$ is a strictly positive bounded weight on $(-\infty, \infty)$.  

\medskip

The following result is known for specialists. We include it for completeness. 

\medskip

\begin{Prop}\label{p2}
Let $\{\alpha_k\}$ be a sequence $(-1,0]$. Then the measure $\mu$ generated by $\{\alpha_k\}$ is in the Steklov class $\Ss$ if and only if $\sum|\alpha_k| < \infty$.
\end{Prop}
\beginpf If $\sum|\alpha_k| < \infty$, the Baxter theorem applies and yields $\mu \in \Ss$.  Conversely, assume that $\mu$ is a measure on~$\T$ whose recurrence coefficients $\alpha_k$ lie in $(-1,0]$. By Lemma~\ref{l7} and Lemma \ref{l3} for $j = 0$, we have
\begin{equation}\label{eq14}
\sup_{n\ge 0} \frac{1}{n}\sum_{k=0}^{n-1}h(k) < \infty.
\end{equation}   
Since $\{\alpha_k\} \subset(-1, 0]$, the sequence $\{h(n)\}$ is increasing and we have $h(n) \ge 1$ for all $n \ge 0$. Hence relation \eqref{eq14} is equivalent to $\sup_{n\ge 0} h(n) < \infty$. As at the end of the proof of Theorem \ref{t1}, we have $\sup_{n\ge 0} |2 s_n + \log h(n)| < \infty$. This implies the Baxter condition $\sum_{k \ge 0}|\alpha_k| < \infty$. \qed

\medskip 

\noindent {\bf Example.} Proposition \ref{p2} shows that the sequence $\{-\frac{1}{4(k+2)}\}_{k \ge 0}$ correspond to a measure $\mu\notin\Ss$. On the other hand, this sequence satisfies Szeg\"o condition $\sum|\alpha_k|^2<\infty$, Steklov bound $\inf s_n \ge -\frac{1}{2}\log n + c$, and the discrete Muckenhoupt condition \eqref{eq55}. Let us construct another sequence $\{\alpha_k\}$ satisfying the Szeg\"o condition and the Steklov bound for which condition $\eqref{eq55}$ is violated. To do this, fix $\delta>\frac{1}{2}$ and construct a sequence of disjoint intervals $I_n = [l_n-\frac{1}{2} , r_n + \frac{1}{2}]$ such that $r_n>l_n\ge 1$ are integer numbers, $r_n-l_n$ is comparable to $l_n^\delta \log l_n$, and the number $\kappa_n = r_n-l_n+1$ of integer points in $I_n$ is a multiple of $4$. Take an interval $I_n$ and divide it into four equal intervals $I_{n,1}, \ldots I_{n,4}$, each containing $\kappa_n/4$ integer points. We enumerate these intervals so that $x\le y$ for all $x \in I_{n,i_1}$, $y\in I_{n, i_2}$, $i_1<i_2$. For $k$ in the left subinterval $I_{n,1}$ of $I_n$ set $\alpha_k = \frac{1}{k^\delta}$. For $k\in I_{n,2}$ we define $\alpha_k$ so that the resulting sequence in $I_{n, 1} \cup I_{n,2}$ is odd with respect to the common (half-integer) point of $I_{n,1}$ and $I_{n,2}$. Then define $\alpha_k$ on $I_{n,3} \cup I_{n,4}$ to obtain the even sequence on $I_n$ with respect to the center on $I_n$. Note that 
$$
\sum_{k \in I_{n,1} \cup I_{n,2}}\alpha_k = 0, \qquad \sum_{k \in I_{n,3} \cup I_{n,4}}\alpha_k = 0, \qquad \sum_{k \in I_n}\alpha_k = 0, \qquad n \ge 1.
$$ 
Finally, for $k$ not in the union $\cup_{n\ge 1}I_n$ we set $\alpha_k = 0$. By construction, the sequence $s_k = \alpha_0 + \ldots + \alpha_k$ is odd on each interval $I_{n}$ with respect to its center and $s_{l_n-1} =0$, $s_{r_n} = 0$ for all $n \ge 1$. In particular, we have
$$
\langle s\rangle_{I_n} = \frac{1}{|I_n|}\sum_{k \in I_n} s_k = 0
$$ 
for the means of $s = \{s_k\}$. Due to symmetricity, we also have
$$
\frac{1}{|I_n|}\sum_{k \in I_n} |s_k| = \frac{1}{|I_{n,1}|}\sum_{k \in I_{n,1}}s_k = 
\frac{1}{|I_{n,1}|}\sum_{k \in I_{n,1}}\sum_{l_n \le s \le k} \frac{1}{s^\delta},
$$
which is comparable to
\begin{equation}\label{eq15}
\frac{1}{r_n-l_n}\int_{l_n}^{r_n}\int_{l_n}^{k}\frac{1}{x^\delta}\,dx = \frac{1}{(1-\delta)}\left(\frac{r_{n}^{2-\delta} - l_{n}^{2-\delta}}{(r_n-l_n)(2-\delta)} - l_{n}^{1-\delta}\right)
\end{equation}
with constants depending only on $\delta$. Our choice of $l_n$ and $r_n$ allows us to estimate the right hand side of \eqref{eq15} from below by $c_1(\delta) \log l_n$, where $c_1(\delta)$ is a constant depending on~$\delta$. Thus, we see that 
\begin{equation}\label{eq16}
\sup_{n\ge 1}\frac{1}{|I_n|}\sum_{k \in I_n} |s_k - \langle s\rangle_{I_n}| = 
\sup_{n\ge 1}\frac{1}{|I_n|}\sum_{k \in I_n} |s_k| \ge 
c_1(\delta)\sup_{n\ge 1}\frac{1}{|I_n|}\log l_n = +\infty.
\end{equation}  
Using Jensen's inequality and the fact that $e^{|x|} \le e^x + e^{-x}$, we obtain 
$$
e^{\frac{2}{|I_n|}\sum_{k \in I_n} |s_k - \langle s\rangle_{I_n}|} 
\le 2 \left(\frac{1}{|I_n|}\sum\nolimits_{k \in I_n} e^{2s_k}\right) \left(\frac{1}{|I_n|}\sum\nolimits_{k \in I_n}e^{-2s_k}\right).
$$
From Theorem~\ref{t1} and \eqref{eq16} we see that one of the measures $\mu$, $\mu_{-1}$ generated by $\{\alpha_k\}$, $\{-\alpha_k\}$, correspondingly, is not in the Steklov class $\Ss$. On the other hand, we have $|\alpha_k| \le \frac{c_2(\delta)}{k^\delta}$ for all $k \ge 1$ and a constant $c_2(\delta)$, hence $\sum|\alpha_k|^2<\infty$. We also have 
$$
|s_k| \le \sum_{k \in I_{n,1}}\frac{1}{k^\delta} \le c_3(\delta)\log k
$$
for all $k\in I_n$ and $s_k = 0$ for $k \notin \cup_{n\ge 1} I_n$. Multiplying, if needed, the elements of $\{\alpha_k\}$ by a small constant, we can obtain a new sequence for which the Steklov bound $s_k \ge -\frac{1}{2}\log k + c$ is satisfied. This ends the construction. 
\bibliographystyle{plain} 
\bibliography{bibfile}
\enddocument